\documentclass{amsart}
\usepackage{amsfonts}

\newtheorem{theorem}{Theorem}
\theoremstyle{plain}

\newtheorem{lemma}{Lemma}

\newtheorem{remark}{Remark}

\numberwithin{equation}{section}

\begin{document}
\title[Integral representation and computation a multiple sum]
{Integral representation and computation a multiple sum in the theory of
cubature formulas}
\author{Georgy\ P.\ Egorychev}
\address{RUSSIA,\ Krasnoyarsk,}
\email{anott@scn.ru}
\date{December 3, 2011}

\begin{abstract}
In this article we furnish a new simple proof of a hard identity from the theory
of cubature formulas via the method of coefficients.
\end{abstract}

\maketitle

\section{{Introduction}}

Let $\mathbf{\alpha }=(\alpha _{0},\alpha _{1},\ldots ,\alpha _{d}),$ $%
\mathbf{\beta }=(\beta _{0},\beta _{1},\ldots ,\beta _{d})$ be vectors from
$E^{d+1}$ with integer non-negative coordinates, and the vector $\mathbf{
\gamma }=(\gamma _{0},\gamma _{1},\ldots ,\gamma _{d})\in E^{d+1}$.
Denote
\textbf{\ }$\left\vert \mathbf{\alpha }\right\vert :=\alpha _{0}+\alpha
_{1}+\ldots +\alpha _{d}=2s+1,$ $\mathbf{\alpha }!:=\alpha _{0}!\alpha
_{1}!\ldots \alpha _{d}!,$ $\binom{\mathbf{\alpha }}{\mathbf{\beta }}:=%
\binom{\alpha _{0}}{\beta _{0}}\ldots \binom{\alpha _{d}}{_{d}},$ where $%
\binom{a}{b}:=\frac{\Gamma \left( a+1\right) }{\Gamma \left( b+1\right)
\Gamma \left( a-b+1\right) },$ and $\binom{a}{b}:=0$, if $b\geq a+1.$
Moreover, we write $\mathbf{\alpha }-1/2:=\left( \alpha _{0}-1/2,\alpha
_{1}-1/2,\ldots ,\alpha _{d}-1/2\right) .$

Heo S. and Xu Y.\cite[pp.631-635]{Heo99} with the help of theory of operators
and generating functions have proved the following multiple combinatorial
identity \cite[the identity (2.9)]{Heo99}:%
$$2^{2s}\mathbf{\alpha }!\binom{\mathbf{\alpha +\gamma }}{\mathbf{\alpha }}%
=
$$
\begin{equation}
\sum_{j=0}^{s}\left( -1\right) ^{j}\binom{d\text{ }+\sum_{i=0}^{d}(\alpha
_{i}+\gamma _{i})}{j}\sum_{\beta _{0}+\beta _{1}+...+\beta
_{d}=s-j}\prod_{i=0}^{d}\binom{\beta _{i}\mathbf{+}\gamma _{i}}{%
\beta _{i}}\left( 2\beta _{i}+\gamma _{i}+1\right) ^{\alpha _{i}}.
\label{K1}
\end{equation}

At the end of the 1970's, G.P. Egorychev has developed the method of
coefficients, which was successfully applied to many
combinatorial sums \cite{2,3,4,5}. The purpose of this article
is finding a new simple proof of identity (\ref{K1}) by means of the
method of coefficients \cite{2} and multiple applications of a known
theorem on the total sum of residues in the theory of holomorphic functions.

\section{{Proof of the identity\ }(\protect\ref{K1})}

The identity (\ref{K1}) can be expressed in the form of%
\begin{equation*}
\sum_{j=0}^{s}\left( -1\right) ^{j}\binom{d+\sum_{i=0}^{d}(\alpha
_{i}+\gamma _{i})}{j}\sum_{\beta _{0}+\beta _{1}+...+\beta
_{d}=s-j}\prod_{i=0}^{d}\binom{\beta _{i}\mathbf{+}\gamma _{i}}{%
\beta _{i}}\frac{\left( 2\beta _{i}+\gamma _{i}+1\right) ^{\alpha _{i}}}{%
\left( \alpha _{i}\right) !}=
\end{equation*}%
\begin{equation}
=2^{2s}\prod_{i=0}^{d}\binom{\alpha _{i}\mathbf{+}\gamma _{i}}{%
\alpha _{i}}.  \label{K2}
\end{equation}%
Denote by $T\left( s;\mathbf{\alpha },\mathbf{\beta }\right)$ the left hand side
of identity (\ref{K2}):%
\begin{equation}
T\left( s;\mathbf{\alpha },\mathbf{\beta }\right) :=\sum_{j=0}^{s}\left(
-1\right) ^{j}\binom{\sum_{i=0}^{d}(\alpha _{i}+\gamma _{i})\mathbf{+}d}{j}%
\times S_{j},  \label{K3}
\end{equation}%
where%
\begin{equation}
S_{j}:=\sum_{\beta _{0}+\beta _{1}+...+\beta _{d}=s-j}\prod_{i=0}^{d}%
\binom{\beta _{i}\mathbf{+}\mu _{i}}{\beta _{i}}\frac{\left( 2\beta _{i}+\mu
_{i}+1\right) ^{\alpha _{i}}}{\left( \alpha _{i}\right) !}.  \label{K4}
\end{equation}%
Then by means of the method of coefficients we obtain%
\begin{equation*}
S_{j}=\sum_{\left\vert \mathbf{\beta }\right\vert
=s-j}(\prod_{i=0}^{d}\binom{\beta _{i}\text{ }\mathbf{+}\text{ }%
\gamma _{i}}{\beta _{i}}\frac{\left( 2\beta _{i}+\gamma _{i}+1\right)
^{\alpha _{i}}}{\left( \alpha _{i}\right) !}=
\end{equation*}%
\begin{equation*}
=\sum_{\mathbf{\beta }_{0}=0}^{\infty }\ldots \sum_{\mathbf{\beta }%
_{d}=0}^{\infty }\mbox{\bf res}_{z_{0},\ldots
,z_{d},t}(t^{-s+j-1}\prod_{i=0}^{d}\left( 1-tz_{i}\right) ^{-\gamma
_{i}-1}z_{i}^{-\beta _{i}-1})\times
\end{equation*}
\begin{equation*}
\times \mbox{\bf res}_{w_{0},\ldots
,w_{d}}(\prod_{i=0}^{d}w_{i}^{-\alpha _{i}-1}\exp (w_{i}\left(
2\beta _{i}+\gamma _{i}+1\right) )=
\end{equation*}%
\begin{equation*}
=\mbox{\bf res}_{w_{0},\ldots
,w_{d},t}\{t^{-s+j-1}(\prod_{i=0}^{d}w_{i}^{-\alpha _{i}-1}\exp
(w_{i}\left( \gamma _{i}+1\right) )\times
\end{equation*}
\begin{equation*}
\times \prod_{i=0}^{d}\left(
\sum_{\beta _{i}=0}^{\infty }\left( \exp (\beta _{i}\left( 2w_{i}\right)
\right) \mbox{\bf res}_{z_{i}}\left( \left( 1-tz_{i}\right) ^{-\gamma
_{i}-1}z_{i}^{-\beta _{i}-1}\right) \right) \}=
\end{equation*}%
(the summation by each $\beta _{i}$ and $\mbox{\bf res}_{z_{i}},$ $%
i=0,\ldots ,d$: the substitution rule, the changes $z_{i}=\exp (2w_{i}),$ $i$
$=$ $0,\ldots ,d$)%
\begin{equation*}
=\mbox{\bf res}_{w_{0},\ldots
,w_{d},t}\{t^{-s+j-1}(\prod_{i=0}^{d}w_{i}^{-\alpha _{i}-1}\exp
(w_{i}\left( \gamma _{i}+1\right) )\times \prod_{i=0}^{d}\left(
1-t\exp (2w_{i})\right) ^{-\gamma _{i}-1}\}=
\end{equation*}%
\begin{equation*}
=\mbox{\bf res}_{w_{0},\ldots
,w_{d},t}\{t^{-s+j-1}\prod_{i=0}^{d}w_{i}^{-\alpha _{i}-1}\left(
\exp (-w_{i})-t\exp (w_{i})\right) ^{-\gamma _{i}-1}\},
\end{equation*}%
i.e.%
\begin{equation}
S_{j}=\mbox{\bf res}_{w_{0},\ldots
,w_{d},t}\{t^{-s+j-1}\prod_{i=0}^{d}w_{i}^{-\alpha _{i}-1}\left(
\exp (-w_{i})-t\exp (w_{i})\right) ^{-\gamma _{i}-1}\}.  \label{K5}
\end{equation}%
According to (\ref{K3})\textbf{--}(\ref{K5})\textbf{\ }we obtain%
\begin{equation*}
T\left( s;\mathbf{\alpha },\mathbf{\beta }\right) =\sum_{j=0}^{s}\mbox{\bf
res}_{w_{0},\ldots
,w_{d},t}\{t^{-s+j-1}\prod_{i=0}^{d}w_{i}^{-\alpha _{i}-1}\left(
\exp (-w_{i})-t\exp (w_{i})\right) ^{-\gamma _{i}-1}\}\times
\end{equation*}
\begin{equation*}
\times \mbox{\bf res}%
_{x}\{x^{-j-1}\left( 1-x\right) ^{d+\sum_{i=0}^{d}(\alpha _{i}+\gamma
_{i})}\}=
\end{equation*}%
\begin{equation*}
=\mbox{\bf res}_{w_{0},\ldots
,w_{d},t}\{t^{-s-1}\prod_{i=0}^{d}w_{i}^{-\alpha _{i}-1}\left( \exp
(-w_{i})-t\exp (w_{i})\right) ^{-\gamma _{i}-1}\times
\end{equation*}
\begin{equation*}
\times(\sum_{j=0}^{\infty }t^{j}%
\mbox{\bf res}_{x}\{x^{-j-1}\left( 1-x\right) ^{d+\sum_{i=0}^{d}(\alpha
_{i}+\gamma _{i})})\}=
\end{equation*}%
(summation w.r.t. $j$, and $\mbox{\bf res}_{x}$: the substitution rule, the
change $x=t$)%
\begin{equation*}
=\mbox{\bf res}_{w_{0},\ldots
,w_{d},t}\{t^{-s-1}\prod_{i=0}^{d}w_{i}^{-\alpha _{i}-1}\left( \exp
(-w_{i})-t\exp (w_{i})\right) ^{-\gamma _{i}-1}\left( 1-t\right)
^{d+\sum_{i=0}^{d}(\alpha _{i}+\gamma _{i})}\}.
\end{equation*}%
Thus we proved

\begin{lemma}
\label{AprL1}Let parameters $s,\alpha _{0},\alpha _{1},\ldots ,\alpha
_{d},\beta _{0},\beta _{1},\ldots ,\beta _{d}$ be non-negative integers,
for which $\alpha _{0}+\ldots +\alpha _{d}=2s+1,$ \textit{and the }%
vector $(\mu _{0},\mu _{1},\ldots ,\mu _{d})\in \mathbb{R}^{d+1}.$
\textit{Then the following integral formula holds}:%
\begin{equation*}
\sum_{j=0}^{s}\left( -1\right) ^{j}\binom{d+\sum_{i=0}^{d}(\alpha
_{i}+\gamma _{i})}{j}\sum_{\beta _{0}+\beta _{1}\ldots +\beta
_{d}=s-j}\prod_{i=0}^{d}\binom{\beta _{i}\mathbf{+}\gamma _{i}}{%
\beta _{i}}\frac{\left( 2\beta _{i}+\gamma _{i}+1\right) ^{\alpha _{i}}}{%
\left( \alpha _{i}\right) !}=
\end{equation*}%
\begin{equation}
=\mbox{\bf res}_{w_{0},\ldots
,w_{d},t}\{t^{-s-1}\prod_{i=0}^{d}w_{i}^{-\alpha _{i}-1}\left( \exp
(-w_{i})-t\exp (w_{i})\right) ^{-\gamma _{i}-1}\left( 1-t\right)
^{d+\sum_{i=0}^{d}(\alpha _{i}+\gamma _{i})}\}.  \label{K6}
\end{equation}
\end{lemma}

\begin{remark}
It is easy to see, that a consecutive calculation of multiple integral in
the right hand side of \textbf{(}\ref{K6}\textbf{) }  on each variable $t
$ \textit{and} $w_{0},\ldots ,w_{d}$ gives the multiple sum in the left part
of (\ref{K6}). Now, we provide new proof of identity (\ref{K2}\textbf{)}
by calculation of a multiple residue at zero point in the right
part of the formula (\ref{K6}) consecutively on each variable $w_{0},\ldots
,w_{d}$ \textit{and} $t$ (see lemmas\textit{\ }\ref{AprL1}--\ref{AprL3} and
the theorem\emph{\ }\ref{AprT1}).
\end{remark}

Let's introduce some necessary notations. Denote
\begin{equation}
f=f\left( w,t\right) :=e^{-w}-te^{w},\text{ }g=g\left( w,t\right)
:=e^{-w}+te^{w},  \label{K7}
\end{equation}%
where $\alpha $ is the fixed integer and $\gamma \in \mathbb{R}.$ Obviously%
\begin{equation}
f^{\prime }:=\frac{df}{dw}=-g,\text{ }g^{\prime }:=\frac{dg}{dw}=-f,\text{ }%
g^{2}-f^{2}=4t,\text{ }(f^{-\gamma })^{\prime }=\gamma f^{-\gamma -1},\text{
}g^{\alpha }=-\alpha g^{\alpha -1}f,  \label{K11}
\end{equation}%
\begin{equation}
f\left( 0\right) =1-t,\text{ }g\left( 0\right) =1+t.  \label{K12}
\end{equation}%

\begin{lemma}
\label{AprL2}\textit{If }$s$\textit{\ is a non-negative integer and
}$\gamma \in \mathbb{R},$ in notation \textbf{(}\ref{K7}\textbf{)} \textit{and} \textbf{(}\ref%
{K11}\textbf{)} \textit{the following expansion } of the derivative%
\begin{equation}
(f^{-\gamma -1})_{w}^{\left( \alpha \right) }=\left( \gamma +1\right) \times
\ldots \times \left( \gamma +\alpha \right) f^{-\gamma -\alpha -1}g^{\alpha
}+\sum_{k=1}^{[\alpha /2]}c_{k}\left( \gamma \right) f^{-\gamma +2k-\alpha
-1}g^{\alpha -2k},  \label{K13}
\end{equation}%
\textit{with integer coefficients }$c_{1},c_{2},\ldots ,c_{[\alpha /2]}$
\textit{is valid.} According to (\ref{K12})\ the formula (\ref{K13})\ generates
the following formula%
\begin{equation*}
\mbox{\bf res}_{w}w_{i}^{-\alpha -1}\left( \exp (-w_{i})-t\exp
(w_{i})\right) ^{-\gamma -1}:=[\left( \exp (-w_{i})-t\exp (w_{i})\right)
^{-\gamma -1}]_{w=0}^{\left( \alpha \right) }/\alpha !=
\end{equation*}%
\begin{equation}
=\binom{\alpha +\gamma }{\alpha }\left( 1-t\right) ^{-\gamma -\alpha
-1}\left( 1+t\right) ^{\alpha }(1+\sum_{k=1}^{[\alpha /2]}h_{k}\left( \alpha
,\gamma \right) \left( 1-t\right) ^{2k}\left( 1+t\right) ^{-2k}),
\label{K14}
\end{equation}%
where the rational coefficients $h_{k}\left( \alpha ,\gamma \right)
:=c_{k}\left( \gamma \right) /\alpha !,$ $k=1,\ldots ,[\alpha /2].$
\end{lemma}

\begin{proof}
The formula \textbf{(}\ref{K13}\textbf{)\ } can be easily proved by induction
on parameter $\alpha $. According  to  \textbf{(}\ref{K11}\textbf{)}
we have for initial values $\alpha =1,2,3$:%
\begin{equation*}
(f^{-\gamma -1})^{\prime }=-\left( \gamma +1\right) f^{-\gamma -2}f^{\prime
}=\left( \gamma +1\right) f^{-\gamma -2}g.
\end{equation*}%
\begin{equation*}
(f^{-\gamma -1})^{^{\prime \prime }}=((f^{-\gamma -1})^{\prime })^{\prime
}=(\left( \gamma +1\right) f^{-\gamma -2}g)^{\prime }=\left( \gamma
+1\right) \left( f^{-\gamma -2}\right) ^{\prime }g+\left( \gamma +1\right)
f^{-\gamma -2}f=
\end{equation*}%
\begin{equation*}
=\left( \gamma +1\right) \left( \gamma +2\right) f^{-\gamma -3}g^{2}-\left(
\gamma +1\right) f^{-\gamma -1}.
\end{equation*}%
\begin{equation*}
(f^{-\gamma -1})^{^{\prime \prime \prime }}=((f^{-\gamma -1})^{^{\prime
\prime }})^{\prime }=(\left( \gamma +1\right) \left( \gamma +2\right)
f^{-\gamma -3}g^{2}-\left( \gamma +1\right) f^{-\gamma -1})^{\prime }=
\end{equation*}%
\begin{equation*}
=\left( \gamma +1\right) \left( \gamma +2\right) \left( \gamma +3\right)
f^{-\gamma -4}g^{3}+\left( \gamma +1\right) \left( \gamma +2\right)
f^{-\gamma -3}2gf-\left( \gamma +1\right) ^{2}f^{-\gamma -2}g=
\end{equation*}%
\begin{equation*}
=\left( \gamma +1\right) \left( \gamma +2\right) \left( \gamma +3\right)
f^{-\gamma -4}g^{3}-\left( \gamma +1\right) (3\gamma +2)f^{-\gamma -2}g.
\end{equation*}%
Further, if formula \textbf{(}\ref{K12}\textbf{)} is
valid for current value $\alpha $, with the help of \textbf{(}\ref{K11}\textbf{%
)} we have%
\begin{equation*}
(f^{-\gamma -1})^{^{\left( \alpha +1\right) }}=((f^{-\gamma -1})^{\left(
\alpha \right) })^{\prime }=(\left( \gamma +1\right) \times \ldots \times
\left( \gamma +\alpha \right) f^{-\gamma -\alpha -1}g^{\alpha
}+\sum_{k=1}^{[\alpha /2]}c_{k}\left( \gamma \right) f^{-\gamma +2k-\alpha
-1}g^{\alpha -2k})^{\prime }=
\end{equation*}%
\begin{equation*}
=\left( \gamma +1\right) \times \ldots \times \left( \gamma +\alpha
+1\right) f^{-\gamma -\alpha -2}g^{\alpha +1}-\left( \gamma +1\right) \times
\ldots \times \left( \gamma +\alpha \right) f^{-\gamma -\alpha }\alpha
g^{\alpha -1}+
\end{equation*}%
\begin{equation*}
+\sum_{k=1}^{[\alpha /2]}c_{k}\left( \gamma \right) \left( \gamma -2k+\alpha
+1\right) f^{-\gamma +2k-\alpha -2}g^{\alpha -2k+1}-\sum_{k=1}^{[\alpha
/2]}c_{k}\left( \gamma \right) f^{-\gamma +2k-\alpha }\left( \alpha
-2k\right) g^{\alpha -2k-1}=
\end{equation*}%
(replacement in the first sum of an index $k-1$ by $k)$%
\begin{equation*}
=\left( \gamma +1\right) \ldots \left( \gamma +\alpha +1\right) f^{-\gamma
-\alpha -2}g^{\alpha +1}+(c_{1}\left( \gamma \right) \left( \gamma +\alpha
-1\right) -\alpha \left( \gamma +1\right) \ldots \left( \gamma +\alpha
\right) )f^{-\gamma -\alpha +1}g^{\alpha -1}+
\end{equation*}%
\begin{equation*}
+\sum_{k=0}^{[\alpha /2]-1}c_{k+1}\left( \gamma \right) \left( \gamma
-2k+\alpha -1\right) f^{-\gamma +2k-\alpha }g^{\alpha
-2k-1}-\sum_{k=1}^{[\alpha /2]}c_{k}\left( \gamma \right) f^{-\gamma
+2k-\alpha }\left( \alpha -2k\right) g^{\alpha -2k-1}=
\end{equation*}%
\begin{equation*}
=\left( \gamma +1\right) \ldots \left( \gamma +\alpha +1\right) f^{-\gamma
-\alpha -2}g^{\alpha +1}+(c_{1}\left( \gamma \right) \left( \gamma +\alpha
-1\right) -c_{[\alpha /2]}\times
\end{equation*}%
\begin{equation*}
\times\left( \gamma \right) f^{-\gamma +2[\alpha
/2]-\alpha }\left( \alpha -2[\alpha /2]\right) g^{\alpha -2[\alpha /2]-1}+
\end{equation*}%
\begin{equation*}
+\sum_{k=1}^{[\alpha /2]-1}(c_{k+1}\left( \gamma \right) \left( \gamma
-2k+\alpha -1\right) -\left( \alpha -2k\right) c_{k}\left( \gamma \right)
)f^{-\gamma +2k-\alpha }g^{\alpha -2k-1},
\end{equation*}%
as was to be shown.
\end{proof}

\begin{lemma}
\label{AprL3}If $s$ \textit{is a non-negative integer then the}
following formulas hold:%
\begin{equation}
J=\mbox{\bf res}_{t}\left( 1-t\right) ^{-1}\left( 1+t\right)
^{2s+1}t^{-s-1}=2^{2s},  \label{K16}
\end{equation}%
\begin{equation}
J_{k}=\mbox{\bf res}_{t}\left( 1-t\right) ^{k-1}\left( 1+t\right)
^{2s-k+1}t^{-s-1}=0,\text{ }\forall k=1,\ldots ,2s.  \label{K17}
\end{equation}
\end{lemma}

\begin{proof}
We have%
\begin{equation*}
J=\mbox{\bf res}_{t}\left( 1-t\right) ^{-1}\left( 1+t\right)
^{2s+1}t^{-s-1}:=\mbox{\bf res}_{t=0}\left( 1-t\right) ^{-1}\left(
1+t\right) ^{2s+1}t^{-s-1}=
\end{equation*}%
(the theorem of the full sum of residues)%
\begin{equation*}
=-\mbox{\bf res}_{t=1}\left( 1-t\right) ^{-1}\left( 1+t\right)
^{2s+1}t^{-s-1}-\mbox{\bf res}_{t=\infty }\left( 1-t\right) ^{-1}\left(
1+t\right) ^{2s+1}t^{-s-1}=
\end{equation*}%
(directly by definition of a residue at a corresponding point)%
\begin{equation*}
=[\left( 1+t\right) ^{2s+1}t^{-s-1}]_{t=1}-res_{t=0}\left( 1-1/t\right)
^{-1}\left( 1+1/t\right) ^{2s+1}\left( 1/t\right) ^{-s}(-1/t)^{2}=
\end{equation*}%
\begin{equation*}
=2^{2s+1}-\mbox{\bf res}_{t=0}\left( 1-t\right) ^{-1}\left( 1+t\right)
^{2s+1}t^{-s}=2^{2s+1}-J\Leftrightarrow J=2^{2s+1}-J\Rightarrow J=2^{2s}.
\end{equation*}%
Let $k$ be the any fixed number from set $\{1,\ldots ,2s\}$. Similarly
to the previous case, we have%
\begin{equation*}
J_{k}=\mbox{\bf res}_{t}\left( 1-t\right) ^{k-1}\left( 1+t\right)
^{2s-k+1}t^{-s-1}:=\mbox{\bf res}_{t=0}\left( 1-t\right) ^{k-1}\left(
1+t\right) ^{2s-k+1}t^{-s-1}-
\end{equation*}%
\begin{equation*}
-\mbox{\bf res}_{t=\infty }\left( 1-t\right) ^{k-1}\left( 1+t\right)
^{2s+1}t^{-s-1} =
\end{equation*}%
(as powers of binomials $\left( 1-t\right) ^{k-1}$ and $\left( 1+t\right)
^{2s-k+1}t^{-s-1}$ non-negative at any $k$ from set $\{1,\ldots ,2s\}$)%
\begin{equation*}
 = 0-\mbox{\bf res}_{t=0}\left( 1-1/t\right) ^{k-1}\left( 1+1/t\right)
^{5}\left( 1/t\right) ^{-3}(-1/t)^{2}=
\end{equation*}%
\begin{equation*}
=-\mbox{\bf res}_{t=0}\left( 1-t\right) ^{-s-1}\left( 1+t\right)
^{2s-k+1}t^{-s-1}=-J_{k}\Longleftrightarrow J_{k}=-J_{k}\Longleftrightarrow
J_{k}=0.
\end{equation*}
\end{proof}

\begin{theorem}
\label{AprT1}The identity (\ref{K2})\ \textit{is valid}.
\end{theorem}

\begin{proof}
According to \textbf{(}\ref{K6}\textbf{)} we have the following integral representation
for sum $T\left( s;\mathbf{\alpha },\mathbf{\beta }\right) $ in the left
hand side of identity \textbf{(}\ref{K2}\textbf{):}%
\begin{equation*}
T\left( s;\mathbf{\alpha },\mathbf{\beta }\right) =\mbox{\bf res}%
_{w_{0},\ldots ,w_{d},t}\{t^{-s-1}\left( 1-t\right)
^{d+\sum_{i=0}^{d}(\alpha _{i}+\gamma
_{i})}\prod_{i=0}^{d}w_{i}^{-\alpha _{i}-1}\left( \exp
(-w_{i})-t\exp (w_{i})\right) ^{-\mu _{i}-1}\}=
\end{equation*}%
\begin{equation*}
=\mbox{\bf res}_{t}\{t^{-s-1}\left( 1-t\right) ^{d+\sum_{i=0}^{d}(\alpha
_{i}+\gamma _{i})}(\prod_{i=0}^{d}\mbox{\bf res}_{w_{i}}w_{i}^{-%
\alpha _{i}-1}\left( \exp (-w_{i})-t\exp (w_{i})\right) ^{-\mu _{i}-1})\}
\end{equation*}%
Calculating in the last expression each of residues w.r.t. variables $w_{0},\ldots
,w_{d}$ by the formula \textbf{(}\ref{K14}\textbf{)} we have%
\begin{equation*}
T\left( s;\mathbf{\alpha },\mathbf{\beta }\right) =\mbox{\bf res}%
_{t}\{t^{-s-1}\left( 1-t\right) ^{d+\sum_{i=0}^{d}(\alpha _{i}+\gamma
_{i})}\times
\end{equation*}%
\begin{equation*}
\times (\prod_{i=0}^{d}\binom{\alpha _{i}+\gamma _{i}}{\alpha _{i}}%
\left( 1-t\right) ^{-\gamma _{i}-\alpha _{i}-1}\left( 1+t\right) ^{\alpha
_{i}}(\sum_{k=1}^{[\alpha _{i}/2]}h_{k}\left( \alpha _{i},\gamma _{i}\right)
\left( 1-t\right) ^{2k}\left( 1+t\right) ^{-2k}))\}=
\end{equation*}%
(trivial cancelations under the product sign $\prod_{i=0}^{d}\ldots $
according to the assumption $\sum_{i=0}^{d}\alpha _{i}=2s+1)$)%
\begin{equation}
=\binom{\mathbf{\alpha +\gamma }}{\mathbf{\alpha }}\mbox{\bf res}%
_{t}\{t^{-s-1}\left( 1-t\right) ^{-1}\left( 1+t\right)
^{2s+1}\prod_{i=0}^{d}(1+\sum_{k=1}^{[\alpha _{i}/2]}h_{k}\left(
\alpha _{i},\gamma _{i}\right) \left( 1-t\right) ^{2k}\left( 1+t\right)
^{-2k})\}.  \label{K18}
\end{equation}%
As $[\alpha _{0}/2]+[\alpha _{1}/2]+\ldots +[\alpha _{d}/2]\leq \lbrack
\sum_{i=0}^{d}\alpha _{i}/2]=s,$ it is easy to see, that after
opening the brackets and simplifying the similar terms the product%
\begin{equation*}
\prod_{i=0}^{d}(1+\sum_{k=1}^{[\alpha _{i}/2]}h_{k}\left( \alpha
_{i},\gamma _{i}\right) \left( 1-t\right) ^{2k}\left( 1+t\right) ^{-2k}),
\end{equation*}%
under the sign  $\mbox{\bf res}_{t}$ in (\ref{K18})\textbf{\ }is\textbf{\ }%
is representable in the form of a polynomial%
\begin{equation*}
1+\sum_{k=2}^{2s}\lambda _{k}\left( 1-t\right) ^{k}\left( 1+t\right) ^{-k},
\end{equation*}%
where coefficients $\lambda _{1},\ldots ,\lambda _{2s-1}$ are some fixed
rational numbers. Thus%
\begin{equation*}
T\left( s;\mathbf{\alpha },\mathbf{\beta }\right) =\binom{\mathbf{\alpha
+\gamma }}{\mathbf{\alpha }}\mbox{\bf res}_{t}\{t^{-s-1}\left( 1-t\right)
^{-1}\left( 1+t\right) ^{2s+1}(1+\sum_{k=1}^{2s}\lambda _{k}\left(
1-t\right) ^{k}\left( 1+t\right) ^{-k})\}=
\end{equation*}%
\begin{equation*}
=\binom{\mathbf{\alpha +\gamma }}{\mathbf{\alpha }}\{\mbox{\bf res}%
_{t}\left( 1-t\right) ^{-1}\left( 1+t\right)
^{2s+1}t^{-s-1}+\sum_{k=1}^{2s}\lambda _{k}\mbox{\bf res}_{t}\left(
1-t\right) ^{k-1}\left( 1+t\right) ^{2s-k+1}t^{-s-1}\}=
\end{equation*}%
(calculation of residues in last expression using formulas \textbf{(}\ref%
{K16}\textbf{)} and \textbf{(}\ref{K17}\textbf{)}%
\begin{equation*}
=\binom{\mathbf{\alpha +\gamma }}{\mathbf{\alpha }}\{2^{2s}+\sum_{k=1}^{2s}%
\lambda _{k}\times 0\}=\binom{\mathbf{\alpha +\gamma }}{\mathbf{\alpha }}%
2^{2s}.
\end{equation*}
\end{proof}

\begin{remark}
It would be interesting to know what additional information one can obtain from
the knowledge of the integral representation of the left hand side of
identity (\ref{K1}).
\begin{equation}
J=\mbox{\bf res}_{w_{0},...,w_{d},t}\{t^{-s-1}\prod_{i=0}^{d}w_{i}^{-\alpha _{i}-1}\left( \exp (-w_{i})-t\exp
(w_{i})\right) ^{-\mu _{i}-1}\left( 1-t\right) ^{d+\sum_{i=0}^{d}(\alpha
_{i}+\gamma _{i})}\}/\mathbf{\alpha }!,  \label{K20}
\end{equation}%
For example, the
integral (\ref{K20}) can written in the following form
\begin{equation}
J=\mbox{\bf res}_{t}\{(t^{-s-1}\prod_{i=0}^{d}\mbox{\bf res}%
_{w_{i}}w_{i}^{-\alpha _{i}-1}\left( \exp (-\lambda _{i}w_{i})-t\exp
(\lambda _{i}w_{i})\right) ^{-\gamma _{i}-1}\}/\mathbf{\alpha }!.
\label{K21}
\end{equation}%
The calculation of integral (\ref{K21}) is connected with studying
of the hyperbolic $t$-sine\textit{\ \cite{FoaHun2005}}%
\begin{equation}
\sinh _{t}(x):=(\exp (-x)-t\exp (x))/2,  \label{K22}
\end{equation}%
and the functions $\sinh _{t}^{-\gamma }(x),$ $\gamma \in
\mathbb{N}$, \textit{and}%
\begin{equation}
J_{a,\gamma }\left( t\right) :=\mbox{\bf res}_{z}(z^{-\alpha -1}\left( \exp
(-z)-t\exp (z)\right) ^{-\gamma -1})/\alpha !=\mbox{\bf res}_{z}(z^{-\alpha
-1}\left( \exp (-z)-t\exp (z)\right) ^{-\gamma -1})/\alpha !.  \label{K23}
\end{equation}%
In my opinion, the study of these functions is interesting, including their
combinatorial interpretation and various corresponding relations.
\end{remark}

The author is thankful to E.Zima and I.Kotsireas for useful comments on early drafts
of this paper.

\end{document}